\documentclass[12pt]{article}
\usepackage{amsthm,amsfonts,amssymb,amscd}

\begin{document}

\begin{center} \large \bf Quasi-polynomial valuations of lattice polytopes
\end{center}\vspace{0.5cm}

\centerline{A.V.Pukhlikov}\vspace{0.5cm}

\parshape=1
3cm 11cm \noindent {\small \quad\quad\quad \quad\quad\quad\quad
\quad\quad\quad {\bf }\newline We prove that for any valuation
$\varphi$ of lattice polytopes in ${\mathbb R}^m$, which is
quasi-polynomial with respect to ${\mathbb Z}^m$-translations, and
any lattice polytope $P$, the function $\varphi(nP)$ is
quasi-polynomial in $n$.

Bibliography: 6 items.} \vspace{1cm}


\vspace{1cm}

\begin{flushright} {\it To Askol'd Georgievich Khovanskii}
\end{flushright}

\section*{Introduction}

In this paper we consider valuations
$$
\varphi\colon {\cal P}(\Lambda)\to R
$$
of lattice polytopes in ${\mathbb R}^m$, where $\Lambda={\mathbb
Z}^m$ and ${\cal P}(\Lambda)$ is the set of all convex polytopes
with vertices in $\Lambda$, and $R$ is a field (possibly of
positive characteristic). Following the model of \cite{PukhKh1},
where valuations, polynomial with respect to
$\Lambda$-trans\-lations, were investigated, we study valuations
that are {\it quasi-polynomial} with respect to
$\Lambda$-trans\-lations, and obtain an analog of one of the main
results of \cite{PukhKh1}: that for a quasi-polynomial valuation
$\varphi$ of degree $\leqslant d$ and any $P\in {\cal P}(\Lambda)$
the $R$-valued function $\varphi(nP)$, $n\geqslant 1$, is a
quasi-polynomial of degree $\leqslant d+m$ (Theorem 1.4; see \S 1
for precise definitions). In other words, we extend one of the
main results of \cite{PukhKh1} from polynomial to quasi-polynomial
valuations.

In our work the central object is the ${\mathbb
Q}(z_1,\dots,z_m)$-valued {\it universal valuation}
$$
\omega\colon P\mapsto \sum_{p=(a_1,\dots,a_m)\in P\cap\Lambda}
z_1^{a_1}\dots z_m^{a_m}
$$
mapping $P$ to a Laurent polynomial in $z_*$. This valuation was
introduced and studied in \cite{PukhKh2, Lawrence1991} and in
numerous works since then (for rational polyhedra and usually
under the name of {\it generating series} --- in our case, it is
{\it generating polynomial}); the most recent paper discussing
this huge amount of work and various results on valuations of
polyhedra is \cite{D'AlessandroPapi2026}.

In this paper we never use presentations of lattice polytopes in
terms of cones, but instead apply the old construction of McMullen
\cite{McMullen1977} and Hadwiger \cite{Hadwiger1952}, reducing the
task to $\Lambda$-simplices, for which the universal valuation is
easy to compute explicitly. Any quasi-polynomial valuation can be
expressed in terms of the universal one, which leads to the proof
of Theorem 1.4.

From the technical viewpoint, the special feature of this paper
is, that instead of expressing the indicator ${\bf 1}_P$ of a
convex polytope as a linear combination of indicators of pointed
cones, which leads to power series, we perform all work within the
class of polytopes, in the spirit of \cite{PukhKh1}, making use of
no power series but only (Laurent) polynomials.

The structure of the paper is as follows. In \S 1 we explain
precisely what is meant by quasi-polynomiality and state the main
result. As an application, the number of integral points in the
multiple polytope $nP$, where $P\in {\cal P}({\mathbb Q})$ has
vertices in ${\mathbb Q}^m\subset {\mathbb R}^m$, turns out to be
a quasi-polynomial in $n$. In \S 2 we discuss the classical
construction of McMullen and Hadwiger and refine it, producing a
{\it precise} presentation of the valuation ${\bf 1}_{n\Delta}$
for a $\Lambda$-simplex $\Delta$ as an integral linear combination
of indicators of finitely many $\Lambda$-polytopes translated by
explicitly given sets of lattice points. Finally, in \S 3 we
explicitly compute the universal valuation $\omega$ for a multiple
simplices and apply differential operators
$z_i\frac{\partial}{\partial z_i}$ to obtain quasi-polynomial
valuations and explain what happens when several variables $z_i$
are replaced by the same variable. At the end of this section we
define the {\it zeta-function} of a lattice polytope.


\section{Quasi-polynomials}

{\bf 1.1. Quasi-polynomial maps.} Let $N$ be an abelian group, $R$
a field, $\chi\colon N\to R^*=R\setminus\{0\}$ a homomorphism, and
$L$ an $R$-vector space.

{\bf Definition 1.1.} A map $h\colon N\to L$ is a quasi-polynomial
of degree $\leqslant d$ ($d\in {\mathbb Z}_+$), associated with
$\chi$, if either

(i) $d=0$ and for every $b\in N$ we have $h(a+b)=\chi(b)h(a)$ for
every $a\in N$, or

(ii) $d\geqslant 1$ and for any $b\in N$ the map $h_b\colon N\to
L$,
$$
h_b(a)=h(a+b)-\chi(b)h(a)
$$
is a quasi-polynomial of degree $\leqslant d-1$ in $a$, associated
with $\chi$.

Obviously, for $d=0$ every quasi-polynomial is of the form
$h(a)=\chi(a)c$, where $c=h(0)\in L$ is arbitrary.

Definition 1.1 is a modification of Definition 1 in \cite[\S
2]{PukhKh1} that explained what is a polynomial map $N\to L$ (no
homomorphism $\chi$ involved). Considering the trivial
homomorphism $\chi\colon N\to \{1_R\}$, we make Definition 1.1
identical to \cite[\S 2, Definition 1]{PukhKh1}.

{\bf Proposition 1.2.} {\it A map $h\colon N\to L$ is a
quasi-polynomial of degree $\leqslant d$, associated with $\chi$,
if and only if the map $h^+\colon N\to L$, where
$$
h^+(a)=\chi(a)^{-1}h(a)
$$
is a polynomial of degree $\leqslant d$. If this is the case, then
for any $a_1,\dots, a_k\in N$ there exist some $b_1,\dots, b_l\in
L$ and integer-valued polynomials $p_i(n_1,\dots, n_k)\colon
{\mathbb Z}^k\to {\mathbb Z}$, $1\leqslant i \leqslant l$, of
degree $\leqslant d$, such that}
$$
h\left(\sum_{j=1}^k n_ja_j\right)=\left(\prod_{j=1}^k
\chi(a_j)^{n_j}\right)\left(\sum_{i=1}^l p_i(n_1,\dots,
n_k)b_i\right).
$$

{\bf Proof.} The main claim is shown by induction on $d\geqslant
0$, given that the case $d=0$ is obvious. If $h^+_b(a)$ is a
polynomial of degree $\leqslant d-1$, then explicitly
$$
\chi(a)^{-1}[\chi(b)^{-1}h(a+b)-h(a)]
$$
is a polynomial of degree $\leqslant d-1$ in $a$, so by induction
$$
\chi(b)^{-1}h(a+b)-h(a)=\chi(b)^{-1}[h(a+b)-\chi(b)h(a)]
$$
is a quasi-polynomial of degree $\leqslant d-1$ in $a$, with $b$
fixed, so that $h_b(a)$ is a quasi-polynomial of degree $\leqslant
d-1$ in $a$, too.

Conversely, if $h_b(a)$ is a quasi-polynomial of degree $\leqslant
d-1$, associated with $\chi$, then by induction
$$
\chi(a)^{-1} h_b(a)=\chi(b)[h^+(a+b)-h^+(a)]
$$
is a polynomial of degree $\leqslant d-1$ in $a$. As $b\in N$ is
fixed, we conclude that $h^+_b(a)$ is a polynomial of degree
$\leqslant d-1$, as required. The second claim is \cite[\S 2,
Remark (ii)]{PukhKh1}. Q.E.D. for the proposition.\vspace{0.3cm}


{\bf 1.2. Quasi-polynomial valuations.} Now let $\Lambda\subset
{\mathbb R}^m$ be an additive subgroup of one of the two types:
either $\Lambda=Q^m$ for some subfield $Q\subset {\mathbb R}$
(e.g. $Q={\mathbb Q}$), or $\Lambda={\mathbb Z}^m$ is a discrete
lattice. Below we will consider the subgroups of the second type,
although the theory works for the first one as well. Let ${\cal
P}(\Lambda)$ be the set of all (bounded) convex polytopes with the
vertices in $\Lambda$.

Let $L$ be an $R$-vector space for some field $R$ and $\chi\colon
\Lambda\to R^*$ a homomorphism. As in \cite[\S 1, Definition
2]{PukhKh1}, an $L$-valued valuation on ${\cal P}(\Lambda)$ means
a map $\varphi\colon {\cal P}(\Lambda)\to L$, satisfying the
inclusion-exclusion property.

{\bf Definition 1.3.} A valuation $\varphi$ is {\it
quasi-polynomial} of degree $\leqslant d$ with respect to the
associated homomorphism
$$
\chi\colon \Lambda\to R^*,
$$
if for every $P\in {\cal P}(\Lambda)$ the map
$$
\Lambda\ni\lambda\mapsto \varphi(P+\lambda)\in L
$$
is a quasi-polynomial of degree $\leqslant d$, associated with
$\chi$.

Definition 1.3 is a modification of \cite[\S 2, Definition 1,
B)]{PukhKh1}.

Let $e_1,\dots, e_m$ be the standard basis of ${\mathbb
Z}^m=\Lambda$ (and ${\mathbb R}^m$). Set $\alpha_i=\chi(e_i)$. By
Proposition 1.2, if $\varphi$ is quasi-polynomial of degree
$\leqslant d$ with respect to $\chi$, then for every $P\in {\cal
P}(\Lambda)$ there exist some $b_1,\dots, b_l\in L$ and
integer-valued polynomials $p_i(z_1,\dots, z_m)$, $1\leqslant
i\leqslant l$, of degree $\leqslant d$, such that the value of
$P+\underline{a}$, where $\underline{a}=(a_1,\dots, a_m)\in
{\mathbb Z}^m$, is equal to
$$
\left(\prod_{j=1}^k \alpha_j^{a_j}\right)\left(\sum_{i=1}^l
p_i(a_1,\dots, a_m)b_i\right).
$$
Obviously, if the homomorphism $\chi$ is trivial, that is,
$\alpha_1=\dots=\alpha_m=1_R$, then we get a polynomial valuation
in the sense of \cite{PukhKh1}.

Now we can state the main result of this paper.

{\bf Theorem 1.4.} {\it Let $\varphi\colon {\cal P}(\Lambda) \to
L$ be a quasi-polynomial valuation of degree $\leqslant d$,
associated with a homomorphism $\chi\colon \Lambda= {\mathbb
Z}^m\to R^*$. Then for every polytope $P\in {\cal P}(\Lambda)$
there are elements $\beta_i\in R^*$, $i=1,\dots, M$, elements
$b_{ij}\in L$ and integer-valued polynomials $p_{ij}(n)$, $n\in
{\mathbb Z}$, of degree $\leqslant d+m$, $i=1,\dots, M$,
$j=1,\dots, M_i$, such that
$$
\varphi(nP)=\sum^M_{i=1}\beta_i^n\left(\sum^{M_i}_{j=1} p_{ij}(n)
b_j\right)
$$
for all} $n\in {\mathbb Z}_{\geqslant 1}$.\vspace{0.3cm}


{\bf 1.3. Counting integral points in rational polytopes.} Let us
consider one particular application of the theorem above.

{\bf Example 1.5.} Let $P$ be a rational convex polytope, that is,
a polytope with rational vertices: $P\in {\cal P}({\mathbb Q}^m)$.
Consider the standard valuation
$$
\varphi(P)=\sum_{p\in P\cap {\mathbb Z}^m} 1=\sharp(P\cap {\mathbb
Z}^m),
$$
counting the integral points in $P$. For some $N\geqslant 1$ we
have $P\in {\cal P}({\mathbb Z}^m)$. Let us fix such a value $N$.
Obviously,
$$
\varphi(P)=\sharp(NP\cap N{\mathbb Z}^m).
$$
Let us write down the right hand side as $\psi_N(NP)$. Set
$\zeta=e^{2\pi i/N}$. For an integral vector
$\underline{a}=(a_1,\dots, a_m)\in{\mathbb Z}^m$ the function
$$
f_N(\underline{a})=\frac{1}{N^m}\prod_{\alpha=1}^m
\left(\sum_{j=0}^{N-1} \zeta^{ja_{\alpha}}\right)
$$
is equal to 1 whenever $\underline{a}\in N{\mathbb Z}^m$ and $0$,
otherwise. Therefore,
$$
\varphi(nP)=\psi_N(NnP)=\sum_{\underline{a}\in n(NP)\cap{\mathbb
Z}^m} f_N(\underline{a}).
$$
But $f_N(\underline{a})$ is a sum of homomorphisms
$$
\chi_{N,j_1,\dots,
j_m}=\frac{1}{N^m}\zeta^{(j_1a_1+\dots+j_ma_m)}\colon {\mathbb
Z}^m\to {\mathbb C}^*,
$$
and the valuation
$$
\psi_{N,j_1,\dots, j_m}(Q)=\frac{1}{N^m}\sum_{\underline{a}\in
Q\cap {\mathbb Z}^m}\zeta^{(j_1a_1+\dots+j_ma_m)}
$$
is obviously quasi-polynomial of degree $\leqslant 0$ with respect
to $\chi_{N,j_1,\dots, j_m}$. Therefore, the integer-valued
function $\varphi(nP)$ is a sum of finitely many (complex-valued)
quasi-polynomials of degree $\leqslant m$:
$$
\varphi(nP)=\sum_{i=1}^M \beta_i^n p_i(n)
$$
with $\beta_i\in {\mathbb C}^*$ some powers of $\zeta$ and
$\mathop{\rm deg} p_i\leqslant m$.


\section{A refinement of Hadwiger-McMullen's \\ construction}

{\bf 2.1. The original construction.} We start with a brief
reminder of the construction in \cite{McMullen1977}, that is
essentially based on \cite{Hadwiger1952}, that was used in
\cite{McMullen1977} to obtain results on translation-invariant
valuations, in the first place, the polynomiality of $\varphi(nP)$
for $P\in {\cal P}(\Lambda)$, $n\in{\mathbb Z}_+$, and $\varphi$
translation-invariant.

(1) For an arbitrary polytope $P\in {\cal P}(\Lambda)$,
considering the sum
$$
\sum_{F}(-1)^{\dim P-\dim F}\beta(F,P)\varphi(F)
$$
over all faces $F$ of $P$, including $P$ itself, where
$\beta(F,P)$ is the internal angle of $P$ at $F$, \cite[Section
3]{McMullen1977} reduces the tasks to the case when $\varphi$ is a
{\it simple} valuation, vanishing on all polytopes of dimension
smaller than $m$.

(2) Dissecting a polytope $P$ of full dimension into
full-dimensional simplices $Q_i$, $i\in I$, whose vertices are
some vertices of $P$, and $Q_i\cap Q_j$ has a smaller dimension
for $i\neq j$, \cite{McMullen1977} further reduces the tasks to
the case when $P$ is a (full-dimensional) simplex with the
vertices in $\Lambda$.

(3) Finally, \cite{McMullen1977} gives a presentation (essentially
due to \cite{Hadwiger1952}) of the multiple simplex $n\Delta$,
$n\in {\mathbb Z}_{>0}$, as a union of translates of finitely many
polytopes that are Minkowski sums of several faces of $\Delta$,
with ${\mathbb R}^m$ being the direct sum of the corresponding
linear subspaces. The sets of points in $\Lambda$ that appear in
these translations are described explicitly; they correspond to
the sets of lattice points in multiple standard simplices.

For the details, see the classical paper \cite{McMullen1977}.

In this paper, we will remove the need for step (1). In out
general setting, the angles $\beta(F,P)$ may not be measured at
all (e.g. when $R$ is a field of positive characteristic) or may
not work for the value space $L$. Therefore, we will need a
refinement of Hadwiger-McMullen's construction: a precise formula
for the multiple simplex $n\Delta$, taking into account polytopes
of smaller dimensions.

This will be done in two steps. In order to get the precise
formula, we first re-trace the original Hadwiger-McMullen's
construction (in a somewhat different way). After that, we obtain
the desired presentation that takes into account polytopes of
non-maximal dimension.\vspace{0.3cm}


{\bf 2.2. Dissecting a multiple simplex.} Now we start using the
language and notations of convex $\Lambda$-chains, introduced in
\cite{PukhKh1}: the Minkowski addition is interpreted as the
convolution of convex chains, which are finite linear combinations
$a_1{\bf 1}_{P_1}+\dots +a_k{\bf 1}_{P_k}$, $a_i\in {\mathbb Z}$
and $P_1,\dots, P_k\in {\cal P}(\Lambda)$, with respect to
integration over the Euler characteristic. The algebra of convex
$\Lambda$-chains is denoted by ${\cal Z}(\Lambda)$ and the
multiplication (generated by the Minkowski addition) is written as
$*$. The subalgebra of zero-dimensional convex chains (where the
polytopes $P_1,\dots, P_k$ are all points) is the group algebra
${\mathbb Z}[\Lambda]$. For a point $\alpha\in \Lambda$ the
indicator ${\bf 1}_{\{\alpha\}}$ will be written as $[\alpha]$.
Obviously, ${\bf 1}_{P}*[\alpha]={\bf 1}_{P+\alpha}$ is the
indicator of the shifted polytope $P+\alpha$.

For every $i\in \{0,\dots, m-1\}$ we denote by ${\cal
Z}_i(\Lambda)$ the additive subgroup of ${\cal Z}(\Lambda)$,
generated by indicators ${\bf 1}_{P}$ of polytopes of dimension
$\leqslant i$. In particular, ${\cal Z}_0(\Lambda)={\mathbb
Z}[\Lambda]$. This is the only case when we get a subalgebra.

Let $x_0, x_1,\dots, x_m$ be a set of points in $\Lambda$ such
that the vectors $x_1-x_0,\dots, x_m-x_0$ are linearly independent
and so form a basis of ${\mathbb R}^m$. Then the simplex with the
these $m+1$ vertices is of full dimension $m$. We will denote this
simplex by $S(x_0,\dots, x_m)$ or $S(\underline{x})$ and its
indicator ${\bf 1}_{S(\underline{x})}$ by $[x_0,\dots,x_m]$ or
just $[\underline{x}]$. For any subset $\{x_{i_1},\dots,
x_{i_k}\}$ of vertices the corresponding simplex will be written
as $S(x_{i_1},\dots, x_{i_k})$ and its indicator as
$[x_{i_1},\dots, x_{i_k}]$.

Hadwiger-McMullen's construction gives a dissection of the
multiple simplex $S(n\underline{x})=S(nx_0,\dots, nx_m)$. We will
reproduce this construction in terms of the indicator
$$
[n\underline{x}]=[\underline{x}]^n=
\underbrace{[\underline{x}]*\cdots*[\underline{x}]}_{n}.
$$
Fix $s\in \{1,\dots, m\}$. For any subset
$$
\{i_0,\dots, i_s\}\subset \{0,1,\dots, m\},
$$
such that
$$
i_0=0\leqslant i_1\leqslant i_2\leqslant\dots\leqslant
i_{s-1}\leqslant  i_s=m-s,
$$
set
$$
\Pi(i_0,\dots, i_s)=[x_0,\dots, x_{i_1}, x_m]*[x_{i_1},\dots,
x_{i_2}, x_{m-1}]*\dots *[x_{i_{s-1}},\dots, x_{i_s}, x_{m+1-s}],
$$
where in the right hand side in every indicator of a simplex
$[x_{i_{a}},\dots, x_{i_{a+1}}, x_{m-a}]$ the vertices $x_j$ run
from $x_{i_{a}}$ to $x_{i_{a+1}}$ consecutively. It is easy to
check that $\Pi(i_0,\dots, i_s)$ is the indicator of a
full-dimensional polytope $P(i_0,\dots, i_s)$, which for
$s\geqslant 2$ is a direct product of $s$ simplices (``cylinders''
in the terminology of \cite{McMullen1977}). When $s=1$, we have
$\Pi(0)=[\underline{x}]$.

Set also for any $s\in \{1,\dots, m\}$
$$
Z(i_0,\dots, i_s)=\sum_{\begin{array}{c} j_1\geqslant 0,\dots,
j_{s+1}\geqslant 0 \\ s+j_1+\dots +j_{s+1}=n+1\end{array}}
[x_m]^{j_1}*\dots *[x_{m-s}]^{j_{s+1}}\in {\mathbb Z}[\Lambda],
$$
where in the right hand side in the product the factors run from
$[x_m]$ to $[x_{m-s}]$ consecutively. In particular,
$$
Z(0)=\sum_{j=1}^n [x_m]^j*[x_{m-1}]^{n-j}.
$$
Now the Hadwiger-McMullen's construction can be presented as
follows.

{\bf Theorem 2.1.} {\it We have}
\begin{equation}\label{28.09.26.1}
[\underline{x}]^{n+1}-\left(\sum_{(i_0,\dots, i_s)}\Pi(i_0,\dots,
i_s)*Z(i_0,\dots, i_s)\right)\in {\cal Z}_{m-1}(\Lambda).
\end{equation}

(This claim is essentially the same as the description of a
dissection of the simplex $S(n\underline{x})$ in
\cite{McMullen1977} up to the difference coming from the
presentation of the standard simplex.)

Since this is a known fact, there is no need to give a proof.
However, we will explain the main steps leading to Theorem 2.1 as
re-tracing them below in Subsection 2.3 will produce a precise
formula for $[\underline{x}]^{n+1}$. First we assume that $x_0=0$
and $e_i=x_i-x_0$ form the standard basis of ${\mathbb R}^m$ with
the corresponding coordinates $t_1,\dots, t_m$. Then we slice the
multiple simplex $S((n+1)\underline{x})$ into {\it storeys}
$$
S_j=S((n+1)\underline{x})\cap\{j\leqslant t_m\leqslant j+1\}.
$$
It is easy to see that $S_0$ can be represented as the union of
the following $m$ polytopes:
$$
\begin{array}{rcl}
S_{0,0} & = & S(x_0,x_m)+S(nx_0,\dots, nx_{m-1}),\\
S_{0,1} & = & S(x_0,x_1,x_m)+S(nx_1,\dots, nx_{m-1}),\\
   & \dots & \\
S_{0,j} & = & S(x_0,\dots, x_j,x_m)+S(nx_j,\dots, nx_{m-1}),\\
   & \dots & \\
S_{0,m-1} & = & S(x_0,\dots,x_{m-1},x_m)+S(nx_{m-1}),
\end{array}
$$
where $S(nx_{m-1})$ is just the point $nx_{m-1}$ and the equality
$$
S_0=S_{0,0}\cup S_{0,1}\cup\dots\cup S_{0,m-1}
$$
is a dissection of $S_0$, that is, the intersection of any two
polytopes in the right hand side is of dimension $\leqslant m-1$.

Doing the same for every storey and combining the results, we get
\begin{equation}\label{28.09.26.2}
\begin{array}{c}
\displaystyle [\underline{x}]^{n+1}-\left([x_m]^n*[\underline{x}]+
\phantom{\sum_{j=0}^{n-1}}\right.\\
\displaystyle +\sum_{j=0}^{n-1}
[x_m]^j*\left.\left(\sum_{i=0}^{m-1}
[x_0,\dots,x_i,x_m]*[x_i,\dots, x_{m-1}]^{n-j}\right)\right)\in
{\cal Z}_{m-1}(\Lambda).
\end{array}
\end{equation}

In particular, $[\underline{x}]$ comes into the huge sum in
brackets with the factor $Z(0)$; this happens when we take $i=m-1$
in the inner brackets. Note that it is easy to check that the
formula (\ref{28.09.26.2}) holds for any $x_0,\dots, x_m\in
\Lambda$, that is, it does not depend on the assumption that
$x_0=0$. Now we have to iterate this formula, applying it to the
factors $[x_i,\dots, x_{m-1}]^{n-j}$. At most $m$ iterations are
needed. The integral parameter $s$ in (\ref{28.09.26.1}) is
precisely the number of iterations until we stop when all powers
are equal to 1. This is a sketch of the proof of Theorem 2.1.


{\bf 2.3. A precise formula for the multiple simplex.} In order to
obtain a precise formula, replacing (\ref{28.09.26.2}), let us
look at the zeroth storey $S_0$ again. Set
$\Theta_j=S_{0,0}\cup\dots\cup S_{0,j}$ and observe that
$$
{\bf 1}_{\Theta(j+1)}={\bf 1}_{\Theta(j)}+{\bf 1}_{S_{0,j+1}}-{\bf
1}_{\Theta(j)\cap S_{0,j+1}},
$$
where $\Theta(j)\cap S_{0,j+1}=S(x_0,\dots,
x_j,x_m)+S(nx_{j+1},\dots,nx_{m-1})$. Therefore, we get a {\it
precise} formula for the zeroth storey
$$
{\bf 1}_{S_0}=\sum_{i=0}^{m-1} [x_0,\dots,x_i,x_m]*[x_i,\dots,
x_{m-1}]^{n}-\sum_{i=1}^{m-1} [x_0,\dots,x_{i-1},x_m]*[x_i,\dots,
x_{m-1}]^{n}.
$$
Obviously, for the intersection $S_0\cap S_1$ of the two adjacent
storeys we have
$$
{\bf 1}_{S_0\cap S_1}=[x_m]*[x_0,\dots, x_{m-1}]^{n}.
$$
If we accept that for $i=0$ the symbol
$$
[x_0,\dots,x_{i-1},x_m]
$$
means just $[x_m]$, we can write down
$$
{\bf 1}_{S_0}-{\bf 1}_{S_0\cap S_1}=\sum_{i=0}^{m-1}
[x_0,\dots,x_i,x_m]*[x_i,\dots, x_{m-1}]^{n}-
$$
$$
-\sum_{i=0}^{m-1} [x_0,\dots,x_{i-1},x_m]*[x_i,\dots,
x_{m-1}]^{n}.
$$
Doing the same for every storey and combining the results, we
replace (\ref{28.09.26.2}) by the following fact.

{\bf Proposition 2.2.} {\it We have}
\begin{equation}\label{29.09.26.1}
\begin{array}{c}
\displaystyle [\underline{x}]^{n+1}=[x_m]^n*[\underline{x}]+
\sum_{j=0}^{n-1} [x_m]^j*\left(\sum_{i=0}^{m-1}
[x_0,\dots,x_i,x_m]*[x_i,\dots, x_{m-1}]^{n-j}-\right. \\
\displaystyle -\left.\sum_{i=0}^{m-1}
[x_0,\dots,x_{i-1},x_m]*[x_i,\dots, x_{m-1}]^{n-j}\right).
\end{array}
\end{equation}

Note that $[\underline{x}]$ comes into the right hand side with
the same factor $Z(0)$ as in (\ref{28.09.26.2}).

Now let us iterate (\ref{29.09.26.1}). For every $s\in\{2,\dots,
m\}$ introduce the tuples
$$
\underline{\varepsilon}=(\varepsilon_1,\dots, \varepsilon_{s-1}),
$$
where $\varepsilon_j\in\{0,1\}$ and define the chains
$\Pi(i_0,\dots, i_s;\varepsilon_1,\dots, \varepsilon_{s-1})$ as
the products
$$
[x_0,\dots,x_{i_1-\varepsilon_1},x_m]*[x_{i_1},\dots,
x_{i_2-\varepsilon_2},x_{m-1}]*\dots * [x_{i_{s-1}},\dots,
x_{i_s},x_{m+1-s}],
$$
where the general factor
$$
[x_{i_a},\dots, x_{i_{a+1}-\varepsilon_{a+1}},x_{m-a}]
$$
is assumed to be just $[x_{m-a}]$ if $i_a=i_{a+1}$ and
$\varepsilon_{a+1}=1$. Finally, let us define the 0-chains
$Z(i_0,\dots, i_s;\varepsilon_1,\dots, \varepsilon_{s-1})$ as
$$
(-1)^{\varepsilon_1+\dots+\varepsilon_{s-1}}
\sum_{\begin{array}{c} j_1\geqslant 0,\dots, j_{s+1}\geqslant 0 \\
s+j_1+\dots +j_{s+1}=n+1\end{array}} [x_m]^{j_1}*\dots *
[x_{m-s}]^{j_{s+1}}.
$$
Define also the 0-chain $Z^*$ as
$$
(-1)^m \sum_{\begin{array}{c} j_1\geqslant 1,\dots, j_m\geqslant 1,
j_{m+1}\geqslant 0 \\
j_1+\dots +j_{m+1}=n-1\end{array}} [x_m]^{j_1}*\dots
[x_{1}]^{j_{m}}* [x_{0}]^{j_{m+1}}
$$
for $n\geqslant m+1$ and set $Z^*=0$ for $n\leqslant m$.

{\bf Theorem 2.3.} {\it We have}
$$
[\underline{x}]^{n+1}=Z^*+\Pi(0)*Z(0)+\sum_{(i_0,\dots,
i_s;\varepsilon_1,\dots, \varepsilon_{s-1}), s\geqslant 2
}\Pi(i_0,\dots, i_s;\varepsilon_1,\dots,
\varepsilon_{s-1})*Z(i_0,\dots, i_s;\varepsilon_1,\dots,
\varepsilon_{s-1}).
$$

The proof is obtained by iterating (\ref{29.09.26.1}). Note that
all chains $\Pi(i_0,\dots, i_s;\varepsilon_1,\dots,
\varepsilon_{s-1})$ are positive-dimensional: they are indicators
of polytopes of dimension $m-\varepsilon_1-\dots-
\varepsilon_{s-1}$.


\section{Valuations of multiple polytopes}

{\bf 3.1. The universal valuation.} In terms of the algebra ${\cal
Z}(\Lambda)$, the claim of Theorem 1.4 is linear in the argument
$P$: if ${\bf 1}_{P}$ is an integral linear combination of
indicators ${\bf 1}_{Q_i}$ $i\in I$, for some $Q_i\in {\cal
P}(\Lambda)$, and Theorem 1.4 holds for all polytopes $Q_i$, then
it holds for $P$. Therefore, we may assume that $P\in {\cal
P}(\Lambda)$ is a simplex (not necessarily of the full dimension
$m$). Next, we apply Theorem 2.3 and see that we have to show the
quasi-polynomiality of $\varphi(Z^*)$ and
$\varphi(\Pi(\underline{i};\underline{\varepsilon})*Z(
\underline{i};\underline{\varepsilon}))$ in $n$ for every tuple
$(\underline{i};\underline{\varepsilon})$. We consider the second
case for some fixed tuple
$(\underline{i};\underline{\varepsilon})$, the other case being
similar. Using the explicit definition of the 0-chain
$Z(\underline{i};\underline{\varepsilon})$ and the last claim of
Proposition 1.2, we further reduce Theorem 1.4 to the following
claim: let $\alpha_1,\dots, \alpha_{s+1} \in R^*$ be any elements
and $p(y_1,\dots, y_{s+1})\in {\mathbb Q}[y_*]$ a polynomial of
(total) degree $\leqslant d$, taking integral values on ${\mathbb
Z}^{s+1}$, then
\begin{equation}\label{05.10.26.1}
\sum_{\begin{array}{c} j_1\geqslant 0,\dots, j_{s+1}\geqslant 0 \\
j_1+\dots +j_{s+1}=n+1-s\end{array}} \alpha_1^{j_1}*\dots *
\alpha_{s+1}^{j_{s+1}}p(j_1,\dots, j_{s+1})
\end{equation}
is quasi-polynomial in $n$ (in the sense of Theorem 1.4). Using
the linearity again, we can further reduce this claim to the case
when
\begin{equation}\label{05.10.26.2}
p(j_1,\dots, j_{s+1})=j_1^{b_1}\dots j_{s+1}^{b_{s+1}},
\end{equation}
$\underline{b}=(b_1,\dots, b_{s+1})\in {\mathbb Z}_+^{s+1}$, is
just a monomial. The quasi-polynomial presentation of
(\ref{05.10.26.1}) depends considerably of whether the values
$\alpha_1,\dots, \alpha_{s+1}$ are distinct or not. We will study
the case of ``general position'' when $\alpha_1,\dots,
\alpha_{s+1}$ are all distinct and then explain what happens for
an arbitrary set $\alpha_1,\dots, \alpha_{s+1}$.

After these reductions we introduce the universal valuation (or
generating polyno\-mial). Recall that $\Lambda={\mathbb Z}^m$ is
the standard lattice in ${\mathbb R}^m$. We say that the map
$\omega\colon {\cal P}(\Lambda)\to {\cal L}(z_1,\dots, z_m)\subset
{\mathbb Q}(z_*)$, where ${\cal L}(z_*)$ is the algebra of Laurent
polynomials in $z_*$ with integral coefficients,
$$
\omega\colon P\mapsto \sum_{\underline{a}\in P\cap\Lambda}
z_1^{a_1}\dots z_m^{a_m},
$$
is the universal quasi-polynomial valuation of degree $\leqslant
0$, associated with the homomorphism
$$
\chi\colon {\mathbb Z}^m \to {\cal L}(z_*)^*
$$
that sends $(a_1,\dots,a_m)$ to $z_1^{a_1}\dots z_m^{a_m}$.

The valuation $\omega$ is indeed quasi-polynomial, as for any
integral vector $\lambda=(\lambda_1,\dots, \lambda_m)\in\Lambda$
$$
\omega(P+\lambda)=z_1^{\lambda_1}\dots z_m^{\lambda_m}\omega(P).
$$
Obviously, the fraction field of ${\cal L}(z_*)$ is the field
${\mathbb Q}(z_1,\dots,z_m)$ of rational functions. In order to
justify the term ``universal'', we show how to obtain an arbitrary
valuation (\ref{05.10.26.1}) with the assumption
(\ref{05.10.26.2}) from the valuation $\omega$. Let us introduce
differential operators
$$
\delta_i=z_i\frac{\partial}{\partial z_i}.
$$
As $\delta_1^{b_1}\dots \delta_m^{b_m}z_1^{a_1}\dots
z_m^{a_m}=a_1^{b_1}\dots a_m^{b_m}z_1^{a_1}\dots z_m^{a_m}$, we
obtain the identity
\begin{equation}\label{12.10.26.1}
\delta_1^{b_1}\dots \delta_m^{b_m}\omega(P)=\sum_{\underline{a}\in
P\cap\Lambda} a_1^{b_1}\dots a_m^{b_m}z_1^{a_1}\dots z_m^{a_m},
\end{equation}
so that indeed every quasi-polynomial valuation (of any degree),
associated with $\chi$, can be obtained from $\omega$, which
justifies its universality.\vspace{0.3cm}


{\bf 3.2. Computations for the standard simplex and applications.}
Let us compute the universal valuation $\omega$ for the standard
simplex
$$
\Delta_m=\{(t_1,\dots, t_m)\in {\mathbb R}_+^m\,|\,
t_1+\dots+t_m=1\}.
$$

{\bf Theorem 3.1.} {\it The following equality holds for
$n\geqslant 1$:}
\begin{equation}\label{07.10.26.1}
\omega(n\Delta_m)=\sum_{i=1}^m\left(\frac{z_i^{m-1}}{\displaystyle
\prod_{j\neq i}(z_i-z_j)}\right) z_i^n.
\end{equation}

This fact is very well known; see, for instance,
\cite{D'AlessandroPapi2026} (with somewhat different notations);
but it is easy to prove by induction on $m\geqslant 1$, the cases
$m=1,2$ being obvious.

The affine version of (\ref{07.10.26.1}) is given by the following
claim.

{\bf Corollary 3.2.} {\it The following equality holds for
$n\geqslant 1$:}
$$
\sum_{\begin{array}{c} (a_1,\dots,a_m)\in {\mathbb Z}_+^m \\
a_1+\dots+a_m\leqslant n \end{array}}z_1^{a_1}\dots z_m^{a_m}=
\sum_{i=1}^m \left(\frac{z_i^{m}}{\displaystyle
(z_i-1)\prod_{j\neq i}(z_i-z_j)}\right)
z_i^n+\frac{1}{(1-z_1)\dots (1-z_m)}.
$$

For $\underline{b}=(b_1,\dots, b_m)\in{\mathbb Z}_+^m$ set
$\delta^{\underline{b}}=\delta_1^{b_1}\dots \delta_m^{b_m}$.
Applying the operator $\delta^{\underline{b}}$ to
(\ref{07.10.26.1}), we compute the valuation
$\delta^{\underline{b}}\omega$ for the multiple simplex
$n\Delta_m$. The result is given below.

{\bf Theorem 3.3.} {\it The following equality holds:
$$
\delta^{\underline{b}}\omega(n\Delta_m)=\sum_{i=1}^m
\sum_{e_i=0}^{b_i} {b_i\choose e_i} R_i(\underline{b},e_i) z_i^n
n^{b_i-e_i},
$$
where
$$
R_i(\underline{b},e_i)=\delta_i^{e_i}\left(\prod_{j\neq
i}\delta_j^{b_j}\right) \left(\frac{z_i^{m-1}}{\displaystyle
\prod_{j\neq i}(z_i-z_j)}\right)\in {\mathbb Q}(z_*)
$$
is a rational function; the highest power of $(z_i-z_j)$ in the
denominator of $R_i(\underline{b},e_i)$ is $1+b_j+e_i$.}

As an application of Theorem 3.1, let us consider an arbitrary
lattice polytope $P$. Denote by ${\cal V}(P)$ the set of its
vertices.

{\bf Theorem 3.4.} {\it There is a set of rational functions
$q_p(z_*)\in {\mathbb Q}(z_*)$, $p\in {\cal V}(P)$, such that for}
$n\geqslant 1$
$$
\omega(nP)=\sum_{p\in {\cal V}(P)} q_p(z_*) z^{np}.
$$

Here for $p=(a_1,\dots,a_m)\in {\mathbb Z}^m$ the symbol
$z^{np}$ means $(z_1^{a_1}\dots z_m^{a_m})^n$.

{\bf Proof.} Representing ${\bf 1}_{P}$ as a ${\mathbb Z}$-linear
combination of indicators of integral simplices, the vertices of
which are in ${\cal V}(P)$, we reduce the claim to the case when
$P$ is a simplex itself. The universal valuation $\omega$ is
associated with $\chi\colon \underline{a}\mapsto
z^{\underline{a}}$, so if
$$
{\bf 1}_{P}=[x_0,x_1,\dots, x_e]
$$
in terms of \S 2, then the values $\alpha_j=z^{x_j}$ are all
distinct. Applying Theorem 3.1 to every term in the presentation
of Theorem 2.3, we obtain the claim of our theorem. Q.E.D.

{\bf Corollary 3.5.} {\it For every $\underline{b}=(b_1,\dots,
b_m)\in{\mathbb Z}_+^m$ we get the presentation
$$
\delta^{\underline{b}}\omega(nP)=\sum_{p\in {\cal V}(P)}
q_{p,\underline{b}}(z_*,n) z^{np},
$$
where for a particular vertex $p=(a_1,\dots,a_m)$ we have
$$
q_{p,\underline{b}}(z_*,n)=\sum_{e_1=0}^{b_1}\dots
\sum_{e_m=0}^{b_m}\left(\prod_{i=1}^m {b_i\choose e_i}
a_i^{b_i-e_i}\right) n^{|\underline{b}|-|\underline{e}|}
\delta^{\underline{e}} (q_p(z_*)),
$$
$|\underline{b}|=b_1+\dots+b_m$, $\underline{e}=(e_1,\dots, e_m)$,
$|\underline{e}|=e_1+\dots +e_m$, so that
$\delta^{\underline{b}}\omega(nP)$ is a quasi-polynomial in $n$ of
degree $\leqslant |\underline{b}|$.}\vspace{0.3cm}


{\bf 3.3. The case of non-general position.} So far we considered
the case of ``general position'' when, in terms of
(\ref{05.10.26.1}), all values are distinct. Let us find out what
happens when some of them are equal. In terms of the formula
(\ref{12.10.26.1}), we formalize this case as follows. Take $k\in
\{1,\dots, m\}$ and any surjective map
$$
\pi\colon\{1,\dots, m\} \to \{1,\dots, k\}.
$$
Let $w_1,\dots, w_k$ be a new set of variables and define a ring
homomorphism
$$
\Phi\colon {\cal L} (z_1,\dots,z_m)\to {\cal L} (w_1,\dots,w_k)
$$
by setting $\Phi (z_i)=w_{\pi(i)}$. For any vector
$\underline{b}=(b_1,\dots, b_m)\in{\mathbb Z}_+^m$ and any $P\in
{\cal P} (\Lambda)$ the image
$\Phi(\delta^{\underline{b}}\omega(nP))\in {\cal L} (w_*)$ is an
integral Laurent polynomial. Let us compute
$\Phi(\delta^{\underline{b}}\omega(nP))$ as a quasi-polynomial in
$n\geqslant 1$ with coefficients in the field ${\mathbb Q}(w_*)$.
As in the case of general position above, we start with the
standard simplex $\Delta_m$.

{\bf Theorem 3.6.} {\it There exist a set of integer-valued
polynomials $p_{ij}(n)$, $i=1,\dots, k$, $j=1,\dots, K(i)$, of
degree}
$$
\mathop{\rm deg} p_{ij}\leqslant
|\pi^{-1}(i)|-1+\sum_{\alpha\in\pi^{-1}(i)}b_{\alpha},
$$
{\it and rational functions $R_{ij}(w_*)$ of the form}
$$
\frac{w_i^{M(i,j)}}{\displaystyle \prod_{j\neq i}
(w_i-w_j)^{c_{ij}}},
$$
{\it such that $M(i,j)\geqslant m-1$ and}
$$
\Phi(\delta^{\underline{b}}\omega(n\Delta_m))=\sum_{i=1}^k w_i^n
\left(\sum_{j=1}^{K(i)} p_{ij}(n) R_{ij}(w_*)\right).
$$

{\bf Proof.} To begin with, we provide a more explicit description
of the rational functions $R_i(\underline{b},e_i)\in {\mathbb
Q}(z_*)$. Using the obvious identity
$$
\delta_j\frac{1}{(z_i-z_j)^{\alpha}}=-\frac{\alpha+1}{(z_i-
z_j)^{\alpha}}+z_i\frac{\alpha+1}{(z_i- z_j)^{\alpha+1}}
$$
for any $\alpha\geqslant 1$ and $j\neq i$, we see that
$$
\left(\prod_{j\neq i}\delta_j^{b_j}\right)
\left(\frac{z_i^{m-1}}{\displaystyle \prod_{j\neq
i}(z_i-z_j)}\right)
$$
(see Theorem 3.3) is an integral linear combination of rational
functions of the form
$$
\frac{z_i^{m-1+|\underline{k}|-k_i}}{\displaystyle \prod_{j\neq
i}(z_i-z_j)^{1+k_j}},
$$
where $\underline{k}=(k_1,\dots, k_m)$, $0\leqslant k_j\leqslant
b_j$ and $|\underline{k}|=k_1+\dots +k_m$. Applying the operator
$\delta_j^{e_j}$ to each of these functions, we obtain a
presentation of $R_i(\underline{b},e_i)$ as an integral linear
combination of rational functions of the form
\begin{equation}\label{19.10.26.1}
z_i^{m-1+|\underline{k}|-k_j+|\underline{c}_i|} \prod_{j\neq i}
\frac{1}{(z_i-z_j)^{1+k_j+c_{ij}}},
\end{equation}
where $\underline{c}_i=(c_{i1},\dots, c_{im})\in {\mathbb Z}_+^m$,
$c_{ii}=0$ and
$$
\sum_{j=1}^m c_{ij}\leqslant e_i.
$$
For convenience set
$$
R_i=\sum_{e_i=0}^{b_i} {b_i\choose e_i} R_i(\underline{b},e_i)
z_i^n n^{b_i-e_i},
$$
$i=1,\dots, k$. Obviously, for $\alpha\not\in \pi^{-1}(i)$ the
value
$$
R_{\alpha}|_{z_j=w_i, \pi(j)=i}
$$
is well defined: we just replace the variables $z_j$, $j\in
\pi^{-1}(i)$, by $w_i$. The ring homomorphism $\Phi\colon {\cal L}
(z_*)\to {\cal L} (w_*)$ can be written as a composition of
elementary homomorphisms $\Phi_{[j\to i]}$ for $j\in \pi^{-1}(i)$
that replace $z_j$ by $w_i$ and do not change other variables. For
instance, if $\pi^{-1}(1)=\{1,2,\dots, l\}$, where $l\geqslant 2$,
then
$$
\Phi_{[1\to 1]}\colon {\cal L} (z_1,z_2,\dots,z_m) \to {\cal L}
(w_1, z_2, \dots,z_m)
$$
replaces $z_1$ by $w_1$ and
$$
\Phi_{[2\to 1]}\colon {\cal L} (w_1,z_2,\dots,z_m) \to {\cal L}
(w_1, z_3, \dots,z_m)
$$
replaces $z_2$ by $w_1$ and does not change $z_3, \dots,z_m$, and
so on. If $k=m$ and $\pi(i)=i$ for all $i\in \{1,\dots, m\}$,
there is nothing to prove.

The operation $\Phi_{[j\to i]}$ can be applied to a rational
function only if the function does not contain $(z_j-w_i)$ in the
denominator. We will consider in detail the case when
$\pi^{-1}(1)=\{1,2,\dots, l\}$, where $l\geqslant 2$, and compute the
image
$$
\Phi_{[l\to 1]}\circ\dots\circ \Phi_{[2\to 1]}\circ \Phi_{[1\to
1]}(\delta^{\underline{b}}\omega(n\Delta_m))\in {\cal L} (w_1,
z_{l+1}, \dots,z_m)
$$
as a rational function in $w_1, z_{l+1} \dots, z_m$. The operation
\begin{equation}\label{20.10.26.1}
\Phi_{[l\to 1]}\circ\dots\circ \Phi_{[1\to 1]}
\end{equation}
is well defined on rational functions in $z_1,\dots,z_m$ that do
not contain $(z_j-z_i)$ for $i,j\in\{1,\dots, l\}$ in the
denominator, so that in the components $R_{l+1},\dots, R_m$ we
just replace the variables $z_1,\dots,z_l$ by $w_1$. Similarly,
after applying (\ref{20.10.26.1}) the subsequent operations
$\Phi_{[j\to i]}$ for $i\geqslant 2$ will be well defined on the
outcome of (\ref{20.10.26.1}) for $R_1+\dots +R_l$.

The operation $\Phi_{[1\to 1]}$ is trivial, just replacing $z_1$
by $w_1$, so that
$$
\Phi_{[1\to 1]}(R_1+\dots +R_l)\in {\mathbb Q}(w_1, z_2,\dots,
z_m).
$$
Set $z_2=w_1+u$. After this substitution we are in the field
${\mathbb Q}(w_1, u, z_3,\dots,z_m)$. As
$$
R_1+R_2+\dots +R_m
$$
is in fact a (Laurent) polynomial, the function
$$
(R_1+R_2)|_{u=0}\in {\mathbb Q}(w_1, z_3,\dots, z_m)
$$
is well defined (for $j\geqslant 3$ the rational functions $R_j$
do not contain $u$ as a factor in the denominator). From the
explicit description of the functions $R_i(\underline{b},e_i)$,
given above, we see that $R_1$ is an integral sum of functions of
the form
$$
\frac{(\mbox{a rational function in}\,\, w_1, z_3,\dots, z_m)}{
u^{1+\alpha}},
$$
with $\alpha\in {\mathbb Z}_+$. Using the presentations
$$
\frac{1}{w_1+u-z_j}=\sum_{\alpha=0}^{\infty}(-1)^{\alpha}\frac{1}{
(w_1-z_j)^{1+\alpha}}u^{\alpha}
$$
and
$$
(w_1+u)^{n+\beta}=\sum_{\gamma=0}^{n+\beta} {n+\beta \choose
\gamma}w_1^{n+\beta-\gamma} u^{\gamma},
$$
we see that in ${\mathbb Q}(w_1,u,z_{\geqslant 3})$ the function
$R_2$ is an integral linear combination of rational functions of
the form
\begin{equation}\label{20.10.26.2}
\frac{n^{b_2-e_2}{n+\beta \choose
\gamma}w_1^{n+\beta-\gamma}u^{\gamma+
|\underline{\alpha}_2|}}{u^{1+\varepsilon}\prod_{j\geqslant 3}
(w_1-z_j)^{1+s_j+\alpha_{2j}}},
\end{equation}
where $\underline{\alpha}_2=(\alpha_{23},\dots,\alpha_{2m})\in
{\mathbb Z}_+^{m-2}$, $|\underline{\alpha}_2|=\alpha_{23}+\dots +
\alpha_{2m}$, $\varepsilon=k_1+c_{21}$, see (\ref{19.10.26.1}),
and
$$
\gamma+|\underline{\alpha}_2|\leqslant 1+\varepsilon,
$$
and a rational function $R_{2,+}\in {\mathbb Q}(w_1,u,z_{\geqslant
3})$, where in fact $R_{2,+}\in u\cdot{\mathbb Q}(w_1,z_{\geqslant
3})[[u]]$. Collecting the components with $\gamma+|\underline{\alpha}_2|<
1+\varepsilon$ and $\gamma+|\underline{\alpha}_2|= 1+\varepsilon$,
respectively, we get the presentation
$$
R_2=R_{2,-}+R_{2,0}+R_{2,+}.
$$
As $R_{2,0}+R_{2,+}\in {\mathbb Q}(w_1,z_{\geqslant 3})[[u]]$, we
conclude that
$$
R_1+R_{2,-}\equiv 0.
$$
Therefore,
$$
\Phi_{[2\to 1]}\circ \Phi_{[1\to 1]}(R_1+\dots
+R_m)=R_{2,0}+(R_3+\dots +R_m)|_{z_1=z_2=w_1}.
$$
Now, equating $z_3,\dots, z_l$ to $w_1$ in a similar way, and
applying $\Phi_{[j\to i]}$ for $i\geqslant 2$, we obtain the
required presentation of
$\Phi(\delta^{\underline{b}}\omega(n\Delta_m))$. It remains to
estimate the degree of this function as a polynomial in $n$. We
use (\ref{19.10.26.1}) again. Looking at the last step, related to
$w_1$, that is, when $\Phi_{[l\to 1]}$ is applied, we get
integer-valued polynomials in $n$ of degree
$$
\gamma+b_l-e_l,
$$
where
$$
\gamma+\sum_{j=l+1}^m \alpha_{lj}=\sum_{j=1}^{l-1} (1+k_j+c_{lj}).
$$
Here $0\leqslant e_l\leqslant b_l$, $0\leqslant k_j\leqslant b_j$,
$c_{i1}+\dots +c_{im}\leqslant e_i$, so that
$$
\gamma+\sum_{j=l+1}^m \alpha_{lj}\leqslant (l-1)+\sum_{j=1}^{l-1}
b_j +e_l.
$$
Adding $b_l-e_l$ and setting $\alpha_{lj}=0$ for $j\geqslant l+1$,
we get the desired upper bound for the degree: it is not higher
than
$$
(l-1)+\sum_{j=1}^{l} b_j,
$$
as claimed by Theorem 3.6. Continuing for the variables $w_2$,
\dots, $w_k$ in exactly the same way, we complete the proof.
Q.E.D.

Now, for a given map $\pi$ and the corresponding homomorphism
$$
\Phi\colon {\cal L} (z_*)\to {\cal L} (w_*)
$$
we immediately obtain the following generalization of Theorem 3.4
and Corollary 3.5 for an arbitrary lattice polytope $P$ (as above,
with the set of vertices ${\cal V}(P)$, and any
$\underline{b}=(b_1,\dots, b_m)\in{\mathbb Z}_+^m$).

{\bf Theorem 3.7.} {\it There is a set of polynomials in $n$
$$
q_{p,\underline{b}}(w_*,n)\in {\mathbb Q}(w_*)[n],
$$
the coefficients of which are rational functions in $w_*$, 
of degree $\leqslant |\underline{b}|+m$, such that
$$
\Phi(\delta^{\underline{b}}\omega(nP))=\sum_{p\in {\cal V}(P)}
q_{p,\underline{b}}(w_*,n)\Phi(z^p)^n.
$$
The polynomials $q_{p,\underline{b}}(w_*,n)$ can be represented as
sums
$$
\sum g_{p,\underline{b},i}(n) R_{p,\underline{b},i}(w_*),
$$
where $g_{p,\underline{b},i}(n)\in {\mathbb Q}[n]$ are
integer-valued polynomials in $n$ and
$R_{p,\underline{b},i}(w_*)\in {\mathbb Q}(w_*)$ some rational
functions.}

Now the arguments of Subsection 3.1 complete the proof of Theorem
1.4.

{\bf Remark 3.8.} In this section we worked with valuations
$$
P\mapsto \sum_{(n_1,\dots,n_m)\in P\cap{\mathbb Z}^m}
n_1^{b_1}\dots n_m^{b_m}z_1^{n_1}\dots z_m^{n_m}\in {\mathbb
Z}[z_*,z^{-1}_*],
$$
where $b_1,\dots, b_m\in {\mathbb Z}_+$. For a polytope
$$
P\subset {\mathbb R}^m_{y_1,\dots, y_m}
$$
with vertices in ${\mathbb Z}^m$, such that
$$
P\subset\{y_1\geqslant 1,\dots, y_m\geqslant 1\},
$$
it seems natural to define its {\it zeta-function}
$$
\zeta_P(s_1,\dots,s_m,z_1,\dots,z_m)=\sum_{(n_1,\dots,n_m)\in
P\cap{\mathbb Z}^m} n_1^{-s_1}\dots n_m^{-s_m}z_1^{n_1}\dots
z_m^{n_m},
$$
with $(s_1,\dots,s_m,z_1,\dots,z_m)\in {\mathbb C}^{2m}$. It is a
polynomial in $z_*$, the coefficients of which are entire
functions in $s_*$. Meromorphic zeta-functions of polyhedra (in
particular, cones) contained in $\{y_*\geqslant 1\}$ can be
defined in a similar way.

{\bf Remark 3.9.} It is easy to see that Theorem 3.1 implies the
Erhardt-type theorem for the universal valuation of the standard simplex.
Using the construction of Theorem 2.3 one can obtain the Erhardt-type
theorem for an arbitrary polytope. A detailed proof will be 
given elsewhere.


\begin{flushleft}
Department of Mathematical Sciences,\\
The University of Liverpool
\end{flushleft}

\noindent{\it pukh@liverpool.ac.uk}

\end{document}